\documentclass{amsart}

\usepackage{algorithm2e,amssymb,mathtools,mleftright,smfthm,stmaryrd,tikz-cd}

\binoppenalty\maxdimen\relpenalty\maxdimen

\title{Faster computation of Witt vectors over polynomial rings}
\author{Rubén Muñoz-\relax-Bertrand}
\address{Université Marie et Louis Pasteur, CNRS, LmB (UMR 6623), 25000 Besançon, France.}

\begin{document}

\begin{abstract}
We describe an algorithm which computes the ring laws for Witt vectors of finite length over a polynomial ring with coefficients in a finite field.
This algorithm uses an isomorphism of Illusie in order to compute in an adequate polynomial ring.
We also give an implementation of the algorithm in SageMath, which turns out to be faster that Finotti's algorithm, which was until now the most efficient one for these operations.
\end{abstract}

\maketitle

Let $p$ be a prime number.
The ring of $p$-typical Witt vectors plays an important role in arithmetic geometry.
Indeed, they satisfy a universal property for $\delta$-rings \cite{joyal}, which can be somehow understood as commutative rings endowed with a lift of the Frobenius endomorphism modulo $p$.
As such, many questions in arithmetic can be translated to characteristic zero with Witt vectors.

They appear for instance in $p$-adic cohomology or in $p$-adic Hodge theory.
But they have also found some applications to the theory of error correcting codes and to cryptography.
In these concrete situations, one needs a software which is able compute with Witt vectors quickly.

Up until recently, only Witt's construction, which is almost one century old, was used for the computation of Witt vectors.
But a recent breakthrough by Finotti gave a very efficient algorithm, when the coefficient ring has characteristic $p$, to compute with Witt vectors \cite{computationswithwittvectorsandthegreenbergtransform}.
Even more recently, the work of Dennerlein is a first step towards faster computations of Witt vectors in mixed characteristic \cite{computationalaspectsofmixedcharacteristicwittvectors}.

\bigskip

Recall that when the coefficient ring is $\mathbb F_p$, then the ring of Witt vectors is isomorphic to $\mathbb Z_p$.
In that case, it is much faster to do the computations in $\mathbb Z_p$, and to roll back to the Witt vectors once they are done.

A similar situation happens when the coefficient ring is a polynomial ring over $\mathbb F_q$, where $q$ is a power of $p$.
This is the situation which is studied in this article.
Illusie gave a description of the ring of Witt vectors in that case \cite{complexedederhamwittetcohomologiecristalline}.
His result is that the ring of truncated Witt vectors is isomorphic to a subring of a polynomial ring, where fast computation algorithms are well-known.

We present here an algorithm which is able to compute Illusie's isomorphism back and forth.
In many situations, especially when the coefficient ring only has one indeterminate, it yields a new algorithm for Witt vectors which is much faster than Finotti's algorithm.

We also give another algorithm by Xavier Caruso which is also efficient.
Both of these algorithms have been implemented in SageMath \cite{sagemath}.
We give a precise comparison between all three algorithms in the last section.
We also explain how these algorithms still work when the base ring is an effective reduced commutative ring of characteristic $p$, without giving an implementation.

\section*{Acknowledgments}

I would like to thank Xavier Caruso for patiently and carefully proofreading all the SageMath code related to Witt vectors I sent him, and also for giving the idea of algorithm \ref{phantom}.
I would also like to thank Frédéric Chapoton for helping me through the process of adding Witt vectors to SageMath during SageDays 128.
So thanks also to the organisers of that workshop, and everyone I have discussed there with.

Many thanks also to the generosity of people at YOYN Cowork in Ålesund, Norway, without whom this paper would have never seen the light of day.

The author is under ``Contrat EDGAR-CNRS no 277952 UMR 6623, financé par la région Bourgogne-Franche-Comté.''.
This work has been financed by ANR-21-CE39-0009-BARRACUDA.
It was also partially written during my position at Université Paris-Saclay, UVSQ, CNRS, Laboratoire de mathématiques de Versailles, 78000, Versailles, France.

Further thanks to be added after the referee process.

\section{Witt vectors}

In all this article, $p$ is a prime number.
Also, $n$ and $m$ will always denote positive natural numbers.

Given a commutative ring $R$, we shall denote by $W\mleft(R\mright)$ its ring of Witt vectors, and by $W_n\mleft(R\mright)$ its ring of truncated Witt vectors of length $n$.
For an introduction on Witt vectors, we refer to \cite[IX. §1]{algebrecommutativehuit}.

We recall that the underlying set of $W\mleft(R\mright)$ is the set $R^\mathbb N$ of sequences in $R$.
Similarly to real numbers, in order to store and compute Witt vectors one can work up to a given precision $n$.
The underlying set of $W_n\mleft(R\mright)$ is the set $R^n$ corresponding to the $n$ first coordinates of a Witt vector.

For each $i\in\mleft\llbracket0,n-1\mright\rrbracket\coloneqq\mleft[0,n-1\mright]\cap\mathbb N$, we let:
\begin{equation*}
F_i\mleft(X_0,\ldots,X_i\mright)\coloneqq\sum_{l=0}^ip^l{X_l}^{p^{i-l}}\in\mathbb Z\mleft[X_0,\ldots,X_i\mright]\text.
\end{equation*}

One then constructs inductively on $i\in\mleft\llbracket0,n-1\mright\rrbracket$ polynomials $S_i$ and $P_i$ in $\mathbb Z\mleft[X_0,\ldots,X_i,Y_0,\ldots,Y_i\mright]$ starting with $S_0\coloneqq X_0+Y_0$ and $P_0\coloneqq X_0Y_0$, and then putting for $i>0$:
\begin{align*}
S_i&\coloneqq\frac{F_i\mleft(X_0,\ldots,X_i\mright)+F_i\mleft(Y_0,\ldots,Y_i\mright)-\sum_{l=0}^{i-1}p^l{S_l}^{p^{i-l}}}{p^i}\text,\\
P_i&\coloneqq\frac{F_i\mleft(X_0,\ldots,X_i\mright)\times F_i\mleft(Y_0,\ldots,Y_i\mright)-\sum_{l=0}^{i-1}p^l{S_l}^{p^{i-l}}}{p^i}\text.
\end{align*}

The ring laws on $W_n\mleft(R\mright)$ are then defined as follows, given two truncated Witt vectors $\mleft(r_i\mright)_{i\in\mleft\llbracket0,n-1\mright\rrbracket},\mleft(s_i\mright)_{i\in\mleft\llbracket0,n-1\mright\rrbracket}\in W_n\mleft(R\mright)$:
\begin{align*}
\mleft(r_i\mright)_{i\in\mleft\llbracket0,n-1\mright\rrbracket}+\mleft(s_i\mright)_{i\in\mleft\llbracket0,n-1\mright\rrbracket}&\coloneqq\mleft(S_0\mleft(r_0,s_0\mright),\ldots,S_{n-1}\mleft(r_0,\ldots,r_{n-1},s_0,\ldots,s_{n-1}\mright)\mright)\text,\\
\mleft(r_i\mright)_{i\in\mleft\llbracket0,n-1\mright\rrbracket}\times\mleft(s_i\mright)_{i\in\mleft\llbracket0,n-1\mright\rrbracket}&\coloneqq\mleft(P_0\mleft(r_0,s_0\mright),\ldots,P_{n-1}\mleft(r_0,\ldots,r_{n-1},s_0,\ldots,s_{n-1}\mright)\mright)\text.
\end{align*}

The recursive definition of the polynomials $\mleft(S_i\mright)_{i\in\mleft\llbracket0,n-1\mright\rrbracket}$ and $\mleft(P_i\mright)_{i\in\mleft\llbracket0,n-1\mright\rrbracket}$ yields an algorithm to compute the ring laws in $W_n\mleft(R\mright)$, when $R$ is an effective commutative ring, that is when its ring operations can be performed by algorithm.
However, this algorithm is very costly.
This is due to the fact that one computes powers of multivariate polynomials with integral coefficients, which themselves have a lot of coefficients.
For instance, for $p=31$ the polynomial $S_2$ has $152\,994$ monomials \cite{computationswithwittvectorsandthegreenbergtransform}.

Of course, these polynomials can be precomputed and stored in a database.
However, one still needs to evaluate them; as the polynomials are huge, this can take a lot of time.

Giving the precise complexity of this algorithm is harder than it seems.
One would need to give the complexity of the exponentiations, which can always be roughly bounded using a dense polynomial representation as in \cite{theefficientcalculationofpowersofpolynomials}.
However, here, our polynomials are very sparse \cite[IX. §1 3]{algebrecommutativehuit}, so one can use the techniques of \cite{parallelsparsepolynomialmultiplicationusingheaps} which have been parallelised \cite{onthecomputationofpowersofsparsepolynomials}.
One can also resort to \cite{onthebitcomplexityofsparsepolynomialandseriesmultiplication} because it should be possible to describe the support of these polynomials, that is the exponents of the monomials with non zero coefficients, and one would then get a formula for the complexity.
This explains why the complexity of this algorithm has never been given.

\bigskip

If we assume that $R$ has characteristic $p$, then one only needs to know the reduction of these polynomials modulo $p$.
In that case, this method computes a lot of useless information, and there is a better way.

Indeed, in that situation Finotti has introduced another algorithm \cite{computationswithwittvectorsandthegreenbergtransform}.
It is the only other algorithm computing the ring laws for Witt vectors in characteristic $p$.
It uses clever combinatorial techniques to give another descriptions of these laws.

Again here, the complexity is not given, for similar reasons as above.
However, both algorithms have been implemented in SageMath \cite{sagemath} by Jacob Dennerlein and the running time is much shorter.

For some rings, such as $\mathbb F_p$, it is a bad idea to employ either the naive or Finotti's algorithm.
Indeed, in that case we have an isomorphism of rings $W_n\mleft(\mathbb F_p\mright)\cong\mathbb Z/p^n\mathbb Z$.
It is therefore much more efficient to compute the ring laws in $\mathbb Z/p^n\mathbb Z$.

In what follows, we will focus on the case where $R$ is a polynomial ring with coefficients in a reduced ring of characteristic $p$.
We shall describe another ring isomorphic to the ring of truncated Witt vectors, and explain how it yields another faster algorithm in that setting.

\section{Computing Illusie's isomorphism}

In this section, $k$ denotes a reduced commutative ring of characteristic $p$.

In what follows, $\underline X$ shall be shorthand for $X_1,\ldots,X_m$ and $\underline{\mleft[X\mright]}$ shall be shorthand for $\mleft[X_1\mright],\ldots,\mleft[X_m\mright]$.
For instance, $k\mleft[X_1,\ldots,X_m\mright]=k\mleft[\underline X\mright]$.

Recall that a $\delta$-ring is a commutative ring $R$ endowed with a map $\delta\colon R\to R$ satisfying some conditions \cite[définition 1]{joyal}.
Here, we will only work in the case where $R$ has no $p$-torsion, in which case those conditions are equivalent to the fact that the map $F\colon\begin{array}{rl}R\to&R\\x\mapsto&x^p+p\delta\mleft(x\mright)\end{array}$ is a morphism of rings \cite[remark 2.2]{prismsandprismaticcohomology}.

\begin{enonce}{Example}
On the ring of polynomials $W\mleft(k\mright)\mleft[\underline X\mright]$, the unique prolongation of the Witt vector Frobenius $F\colon W\mleft(k\mright)\to W\mleft(k\mright)$ satisfying $F\mleft(X_i\mright)={X_i}^p$ for every $i\in\mleft\llbracket1,m\mright\rrbracket$ is a $\delta$-ring.
We will always endow this ring with this $F$.
\end{enonce}

We recall that the $\delta$-ring of Witt vectors satisfies the following universal property.
Given a $\delta$-ring $A$ and a morphism of commutative rings $\varphi\colon A\to B$, then there exists a unique morphism of $\delta$-rings making the following diagram commutative \cite[théorème 4]{joyal}:
\begin{equation*}
\begin{tikzcd}
A\arrow[rd,"\varphi"]\arrow[r,dashed]&W\mleft(B\mright)\arrow[d]\\
&B\text.
\end{tikzcd}
\end{equation*}

The vertical arrow is the map which to a Witt vector associates its first coordinate.

\begin{enonce}{Example}
The morphism of rings $W\mleft(k\mright)\mleft[\underline X\mright]\to k\mleft[\underline X\mright]$ given by the quotient modulo $V\mleft(W\mleft(k\mright)\mright)W\mleft(k\mright)\mleft[\underline X\mright]$ yields a morphism of $\delta$-rings $W\mleft(k\mright)\mleft[\underline X\mright]\to W\mleft(k\mleft[\underline X\mright]\mright)$.

Recall that for every $x\in k$, the representative $\mleft[x\mright]$ satisfies $F\mleft(\mleft[x\mright]\mright)=\mleft[x\mright]^p$ \cite[IX. §1 proposition 5]{algebrecommutativehuit}.
In particular, this implies that the unique morphism of commutative $W\mleft(k\mright)$-algebras $W\mleft(k\mright)\mleft[\underline X\mright]\to W\mleft(k\mleft[\underline X\mright]\mright)$ such that $X_i$ maps to $\mleft[X_i\mright]$ is the unique morphism of $\delta$-rings making the adjunction diagram commutative.
\end{enonce}

In \cite[corollaire 2.15]{complexedederhamwittetcohomologiecristalline}, Illusie gives an isomorphism between $W\mleft(k\mleft[X_1,\ldots,X_m\mright]\mright)$ and a ring of formal power series.
The goal of the next proposition is to explain how to compute this isomorphism on truncated Witt vectors.

\begin{enonce}{Proposition}\label{modularequality}
Consider $G\mleft(\underline X\mright)\in W\mleft(k\mright)\mleft[\underline X\mright]$, and denote by $\overline G\mleft(\underline X\mright)$ its image in $k\mleft[\underline X\mright]$.
Then:
\begin{equation*}
\mleft[\overline G^{p^{n-1}}\mleft(\underline X\mright)\mright]\equiv G^{p^{n-1}}\mleft(\underline{\mleft[X\mright]}\mright)\pmod{V^{n}\mleft(W\mleft(k\mleft[\underline X\mright]\mright)\mright)}\text.
\end{equation*}
\end{enonce}

\begin{proof}
We show this by induction on $n\in\mathbb N^*$, the case $n=1$ being clear.

Then, by induction hypothesis there is $\epsilon\in W\mleft(k\mleft[\underline X\mright]\mright)$ such that we get:
\begin{equation*}
\begin{split}
\mleft[\overline G^{p^n}\mleft(\underline X\mright)\mright]&=\mleft[\overline G^{p^{n-1}}\mleft(\underline{\mleft[X\mright]}\mright)\mright]^p\\
&=\mleft(G^{p^{n-1}}\mleft(\underline{\mleft[X\mright]}\mright)+V^n\mleft(\epsilon\mright)\mright)^p\text.
\end{split}
\end{equation*}

The conclusion then follows from the binomial expansion and \cite[IX. §1 proposition 5]{algebrecommutativehuit}.
\end{proof}

\begin{enonce}{Proposition}\label{illusieisomorphism}
Let $\mleft(\overline{G_i}\mleft(\underline X\mright)\mright)_{i\in\mleft\llbracket0,n-1\mright\rrbracket}\in W_n\mleft(k\mleft[\underline X\mright]\mright)$ be a truncated Witt vector.
Choose any family $\mleft(G_i\mleft(\underline X\mright)\mright)_{i\in\mleft\llbracket0,n-1\mright\rrbracket}\in{W_n\mleft(k\mright)\mleft[\underline X\mright]}^n$ of polynomials whose projections modulo $V\mleft(W_n\mleft(k\mright)\mright)W_n\mleft(k\mright)\mleft[\underline X\mright]$ are the coefficients of the Witt vector.

We then have the following equality in $W_n\mleft(k\mleft[\underline X\mright]\mright)$:
\begin{equation*}
F^{n-1}\mleft(\mleft(\overline{G_0}\mleft(\underline X\mright),\ldots,\overline{G_{n-1}}\mleft(\underline X\mright)\mright)\mright)=\sum_{i=0}^{n-1}p^i{G_i}^{p^{n-1-i}}\mleft(\underline{\mleft[X\mright]}\mright)\text.
\end{equation*}
\end{enonce}

\begin{proof}
We have:
\begin{equation*}
\begin{split}
F^{n-1}\mleft(\mleft(\overline{G_0}\mleft(\underline X\mright),\ldots,\overline{G_{n-1}}\mleft(\underline X\mright)\mright)\mright)&=F^{n-1}\mleft(\sum_{i=0}^{n-1}V^i\mleft(\mleft[\overline{G_i}\mleft(\underline X\mright)\mright]\mright)\mright)\\
&=\sum_{i=0}^{n-1}p^iF^{n-1-i}\mleft(\mleft[\overline{G_i}\mleft(\underline X\mright)\mright]\mright)\text.
\end{split}
\end{equation*}

Conclude with proposition \ref{modularequality}.
\end{proof}

\section{The algorithm}

Let us keep the notations from the previous section.
Assume furthermore that $k$ is effective, and that one can compute the $p$-th root of any element in $k$ having one.

Since $k\mleft[\underline X\mright]$ is a reduced ring of characteristic $p$, the endomorphism of rings $F^{n-1}\colon W_n\mleft(k\mleft[\underline X\mright]\mright)\to W_n\mleft(k\mleft[\underline X\mright]\mright)$ is injective.
Combined with the fact that, in proposition \ref{illusieisomorphism}, the expression of its image as the evaluation of a polynomial with coefficients in $W_n\mleft(k\mright)$ does not depend on the choices of the lifts, we have thus defined an injective map:
\begin{equation*}
\widetilde F^{n-1}\colon W_n\mleft(k\mleft[\underline X\mright]\mright)\to W_n\mleft(k\mright)\mleft[\underline X\mright]\text.
\end{equation*}

For each $u\in\mathbb N^n$, let $\operatorname v_{p,m}\mleft(u\mright)=\max\mleft\{k\in\mleft\llbracket0,n-1\mright\rrbracket\mid\forall i\in\mleft\llbracket1,n\mright\rrbracket,\ u_i\in p^k\mathbb N\mright\}$.
The image of $\widetilde F^{n-1}$ is included in the ring consisting of polynomials of the form $\sum_{u\in\mathbb N^n}a_u\prod_{i=1}^n{X_i}^{u_i}\in W_n\mleft(k\mright)\mleft[\underline X\mright]$ with $a_u\in p^{n-1-\operatorname v_{p,m}\mleft(u\mright)}W_n\mleft(k\mright)$.
It is an equality when $k$ is perfect.

As explained above, the computation of the ring laws of Witt vectors can be quite costly.
However, computations in rings of polynomials are well-known and much faster.
So we give here an algorithm to compute Witt vectors.

\begin{enonce}{Algorithm}\label{illusie}\

\begin{algorithm}[H]\LinesNumbered
\KwIn{
$m,n\in\mathbb N^*$

$\mleft(\overline{G_i}\mleft(\underline X\mright)\mright)_{i\in\mleft\llbracket0,n-1\mright\rrbracket}\in W_n\mleft(k\mleft[\underline X\mright]\mright)$

$\mleft(\overline{H_i}\mleft(\underline X\mright)\mright)_{i\in\mleft\llbracket0,n-1\mright\rrbracket}\in W_n\mleft(k\mleft[\underline X\mright]\mright)$
}
\KwOut{
$\mleft(\overline{G_i}\mleft(\underline X\mright)\mright)_{i\in\mleft\llbracket0,n-1\mright\rrbracket}+\mleft(\overline{H_i}\mleft(\underline X\mright)\mright)_{i\in\mleft\llbracket0,n-1\mright\rrbracket}$

\tcp*[f]{(respectively $\mleft(\overline{G_i}\mleft(\underline X\mright)\mright)_{i\in\mleft\llbracket0,n-1\mright\rrbracket}\times\mleft(\overline{H_i}\mleft(\underline X\mright)\mright)_{i\in\mleft\llbracket0,n-1\mright\rrbracket}$)}
}
\For{$i\in\mleft\llbracket0,n-1\mright\rrbracket$}{
Fix $G_i\mleft(\underline X\mright)\in W_n\mleft(k\mright)\mleft[\underline X\mright]$ a lift of $\overline{G_i}$

Fix $H_i\mleft(\underline X\mright)\in W_n\mleft(k\mright)\mleft[\underline X\mright]$ a lift of $\overline{H_i}$
}
$G\coloneqq\sum_{i=0}^{n-1}p^i{G_i}^{p^{n-1-i}}\mleft(\underline{X}\mright)$

$H\coloneqq\sum_{i=0}^{n-1}p^i{H_i}^{p^{n-1-i}}\mleft(\underline{X}\mright)$

$R\coloneqq G+H$\tcp*[f]{(respectively $R\coloneqq G\times H$)}

\For{$i\in\mleft\llbracket0,n-1\mright\rrbracket$}{
Let $\overline{R_i}$ be the inverse image through $\operatorname{Frob}^{n-1-i}$ of the projection $\overline R$

Fix $R_i\in W_n\mleft(k\mright)\mleft[\underline X\mright]$ a lift of $\overline{R_i}$

$R\coloneqq\frac{R-{R_i}^{p^{n-1-i}}}p$
}
\Return{$\mleft(\overline{R_i}\mleft(\underline X\mright)\mright)_{i\in\mleft\llbracket0,n-1\mright\rrbracket}$}
\end{algorithm}
\end{enonce}

\begin{enonce}{Remark}
Algorithm \ref{illusie} can be generalised to rings $R$ having a surjective morphism of rings $k\mleft[\underline X\mright]\to R$ given by an algorithm, and for which there also exists an algorithm computing a preimage in $R$ of any element of $k\mleft[\underline X\mright]$.
Then, computing the sum or the product of truncated Witt vectors in $W_n\mleft(R\mright)$ is done by choosing preimages in $W_n\mleft(k\mleft[\underline X\mright]\mright)$, doing computations with algorithm \ref{illusie}, and then computing the image of the result in $W_n\mleft(R\mright)$.

When $k$ does not have characteristic $p$, it might not be possible to use similar methods.
When $k$ has characteristic at least $p^p$, it is even impossible.
Indeed, we have $F_1\mleft(p,-p^{p-1}\mright)=0$ where $F_1$ is the ghost polynomial.
In particular, this method does not yield an injective morphism.
\end{enonce}

\section{Another approach}

Let us keep the notations of the previous section.

The algorithm given above can be described as follows.
Starting with Witt vectors in $W_n\mleft(k\mleft[\underline X\mright]\mright)$, we lift them to $w,w'$ in the ring $W_n\mleft(W_n\mleft(k\mright)\mleft[\underline X\mright]\mright)$.
Then, we compute the $n-1$-th ghost polynomials $F_{n-1}\mleft(w\mright)$ and $F_{n-1}\mleft(w'\mright)$ of the lifted Witt vectors.
We do our computations with these polynomials, and in the end we roll back to Witt vectors.

But in fact, we could have simply done the following thing: instead of computing only the $n-1$-th ghost polynomial, we can compute the $n$ first ones (recall that indexing starts at $0$).
Then, we do the computations in the ring of $n$-tuples of polynomials, and compute the preimage of the result.
We detail this in algorithm \ref{phantom} below.

What we lose here is having to compute every ghost component.
However, there is no need to compute the preimages of the Frobenius as we did in the algorithm in the previous section.
That is because in $F_i$ the leading monomial in $X_i$ is $p^iX_i$, so we can retrieve this monomial by subtracting powers of the $F_j$ for $j<i$, starting with $F_0=X_0$.

\begin{enonce}{Algorithm}[Xavier Caruso]\label{phantom}\

\begin{algorithm}[H]\LinesNumbered
\KwIn{
$m,n\in\mathbb N^*$

$\mleft(\overline{G_i}\mleft(\underline X\mright)\mright)_{i\in\mleft\llbracket0,n-1\mright\rrbracket}\in W_n\mleft(k\mleft[\underline X\mright]\mright)$

$\mleft(\overline{H_i}\mleft(\underline X\mright)\mright)_{i\in\mleft\llbracket0,n-1\mright\rrbracket}\in W_n\mleft(k\mleft[\underline X\mright]\mright)$
}
\KwOut{
$\mleft(\overline{G_i}\mleft(\underline X\mright)\mright)_{i\in\mleft\llbracket0,n-1\mright\rrbracket}+\mleft(\overline{H_i}\mleft(\underline X\mright)\mright)_{i\in\mleft\llbracket0,n-1\mright\rrbracket}$

\tcp*[f]{(respectively $\mleft(\overline{G_i}\mleft(\underline X\mright)\mright)_{i\in\mleft\llbracket0,n-1\mright\rrbracket}\times\mleft(\overline{H_i}\mleft(\underline X\mright)\mright)_{i\in\mleft\llbracket0,n-1\mright\rrbracket}$)}
}
\For{$i\in\mleft\llbracket0,n-1\mright\rrbracket$}{
Fix $G_i\mleft(\underline X\mright)\in W_n\mleft(k\mright)\mleft[\underline X\mright]$ a lift of $\overline{G_i}$

Fix $H_i\mleft(\underline X\mright)\in W_n\mleft(k\mright)\mleft[\underline X\mright]$ a lift of $\overline{H_i}$
}
\For{$j\in\mleft\llbracket0,n-1\mright\rrbracket$}{
$P_j\coloneqq\sum_{i=0}^{j}p^i{G_i}^{p^{j-i}}\mleft(\underline{X}\mright)$

$Q_j\coloneqq\sum_{i=0}^{j}p^i{H_i}^{p^{j-i}}\mleft(\underline{X}\mright)$

$T_j\coloneqq P_j+Q_j$\tcp*[f]{(respectively $T_j\coloneqq P_j\times Q_j$)}
}
$R_0\coloneqq T_0$

\For{$j\in\mleft\llbracket1,n-1\mright\rrbracket$}{
$R_j\coloneqq\frac{T_j-\sum_{i=0}^{j-1}p^i{R_i}^{p^{j-1-i}}\mleft(\underline{X}\mright)}{p^j}$
}
\Return{the projection $\mleft(\overline{R_i}\mleft(\underline X\mright)\mright)_{i\in\mleft\llbracket0,n-1\mright\rrbracket}$}
\end{algorithm}
\end{enonce}

In some cases, this approach is faster.
We give comparisons in the next section.

In practice, for instance when $k=\mathbb F_p$, after lifting to $\mathbb Z_p$ it is more efficient to do computations in $\mathbb Q_p$.
In SageMath \cite{sagemath}, which uses FLINT \cite{flint}, the floating point approach yields much better results.
We only need to take care of the precision, but here clearly precision $p^n$ is enough for our needs.

The author is thankful to Xavier Caruso for both giving this another algorithm, and pointing out the fact that computations are faster in $\mathbb Q_p$ rather than in $\mathbb Z_p$.

\begin{enonce}{Remark}
One of the advantages of this approach is that it can be generalised to a wide category of rings, even when the characteristic of the base ring is not $p$.

Indeed, if $R$ is a ring of characteristic $p^l$ for some $l\in\mathbb N^*$, and if we assume that there exists a surjective ring morphism $S\to R$ where $S$ is a $p$-torsion free commutative ring, $S/p^{n+l-1}S$ is effective, and such that an algorithm can compute both images and preimages of the projection, then algorithm \ref{phantom} can be adapted.
Assuming $k$ is perfect here we would have $R=k\mleft[\underline X\mright]$, $l=1$ and $S=W\mleft(k\mright)\mleft[\underline X\mright]$ in algorithm \ref{phantom}.
We omit the details.

This setting encompasses for instance the case $R=\mathcal O_K/p\mathcal O_K$, where $K$ is a ramified extension of $\mathbb Q_p$, which is studied in $p$-adic Hodge theory.

Things can be trickier when studying quotient groups of $\mathbb Z/p^2\mathbb Z\mleft[\underline X\mright]$ for instance.
It might be tempting to raise it to a quotient of $\mathbb Z_p\mleft[\underline X\mright]$ defined with the same equations, but it might have $p$-torsion, as can be seen with $\mathbb Z/p^2\mathbb Z\mleft[X\mright]/\mleft\langle pX\mright\rangle$.
\end{enonce}

For more about computations of Witt vectors in mixed characteristic, see \cite{computationalaspectsofmixedcharacteristicwittvectors}.

\section{Comparison between the algorithms}

To conclude this article, we compare the algorithm of Finotti, and algorithms \ref{illusie} (referred in the tables below as \texttt{Illusie}) and \ref{phantom} (referred in the tables below as \texttt{phantom}).

All of the computations here were made with SageMath 10.6.rc0  \cite{sagemath}, running on a 13th Gen Intel Core i7-13850HX processor with Debian GNU/Linux version 12.10.

SageMath in turn uses various external libraries.
Computations in $p$-adics are done with FLINT \cite{flint}, computations with polynomials sometimes use Singular \cite{singular} and computations in finite fields use Givaro \cite{givaro} for fields of small cardinality and PARI/GP \cite{pari} otherwise.

\bigskip

First, for $d,n\in\mathbb N$, $p$ a prime number and $q$ a power of $p$, we compare the running time, in seconds, and maximum memory usage, in bytes, of each algorithms to add (respectively multiply) two randomly chosen $p$-typical Witt vectors in $W_n\mleft(\mathbb F_q\mleft[X\mright]\mright)$ such that the degree of each coefficient is at most $d$.

Note that Finotti's algorithm precomputes some auxiliary polynomials which are used for computations.
Once these are computed, they can be employed for any future computations.
For that reason, we did not take into account the precomputation time and memory of Finotti's algorithm, even though it is not always negligible.
Therefore, the values below can be larger for Finotti's method.

\bigskip

\begin{figure}[ht]
\begin{tabular}{|c|c|c|c|c|c|c|}
\hline
$+$&\multicolumn{3}{|c|}{Time (s)}&\multicolumn{3}{|c|}{Memory (B)}\\
\hline\hline
$d$&\texttt{Finotti}&\texttt{phantom}&\texttt{Illusie}&\texttt{Finotti}&\texttt{phantom}&\texttt{Illusie}\\
\hline\hline
0&3&$<$1&$<$1&115 988&500 555&40 742\\
10&464&21&21&255 674&1 181 842&357 153\\
20&644&40&57&411 826&1 874 590&689 184\\
30&1 023&143&120&605 842&2 567 053&976 949\\
40&1 517&276&247&718 202&3 249 853&1 343 862\\
50&1 861&370&319&836 682&3 930 183&1 607 118\\
60&2 286&556&445&1 114 602&4 636 076&1 924 099\\
70&2 778&740&597&1 220 330&5 317 147&2 423 381\\
80&3 089&1 109&909&1 339 018&5 994 007&2 663 495\\
90&3 227&1 003&821&1 457 706&6 681 261&2 899 909\\
100&3 959&1 561&1 365&1 576 786&7 381 732&3 178 384\\
\hline
\hline
$\times$&\multicolumn{3}{|c|}{Time (s)}&\multicolumn{3}{|c|}{Memory (B)}\\
\hline\hline
$d$&\texttt{Finotti}&\texttt{phantom}&\texttt{Illusie}&\texttt{Finotti}&\texttt{phantom}&\texttt{Illusie}\\
\hline\hline
0&1&$<$1&$<$1&223 584&495 027&32 117\\
10&129&40&36&447 501&1 860 042&599 806\\
20&238&152&118&689 659&3 217 722&1 173 548\\
30&397&345&263&821 377&4 588 223&1 754 689\\
40&536&736&528&1 182 609&5 957 590&2 323 490\\
50&669&1 136&814&1 284 700&7 328 680&2 899 679\\
60&802&1 640&819&1 442 013&8 705 957&3 484 596\\
70&761&1 958&1 234&1 601 623&10 079 024&4 053 415\\
80&1103&2 910&2 072&2 166 266&11 422 718&4 627 325\\
90&975&2 708&1 835&2 237 539&12 814 396&5 199 606\\
100&1372&2 970&1 276&2 370 325&14 175 216&5 781 967\\
\hline
\end{tabular}
\caption{Comparisons in $W_5\mleft(\mathbb F_{3^{10}}\mleft[X\mright]\mright)$.}\label{d}
\end{figure}

Figure \ref{d} describes the behaviour of the three algorithms for varying $d$.
Even though the computation time for all algorithms is roughly linear in $d$, there is a lot of noisy variations.

For instance, one can notice that algorithm \ref{illusie} always outperforms the two other algorithms for the addition, except for $d=20$.
Similarly, for the multiplication Finotti's algorithm and algorithm \ref{illusie} run in similar timings, except for small $d$ where algorithm \ref{illusie} is faster, and for $50\leqslant d\leqslant90$ where Finotti's algorithm is much faster.

In almost all cases, Finotti's method uses less memory.
This will also be the case in most comparisons below.

\begin{figure}[ht]
\begin{tabular}{|c|c|c|c|c|c|c|}
\hline
$+$&\multicolumn{3}{|c|}{Time (s)}&\multicolumn{3}{|c|}{Memory (B)}\\
\hline\hline
$n$&\texttt{Finotti}&\texttt{phantom}&\texttt{Illusie}&\texttt{Finotti}&\texttt{phantom}&\texttt{Illusie}\\
\hline\hline
1&$<$1&$<$1&$<$1&26 092&522 697&42 155\\
2&$<$1&$<$1&$<$1&99 010&509 745&41 122\\
3&$<$1&$<$1&$<$1&99 114&515 745&57 035\\
4&32&1&1&141 514&655 715&133 229\\
5&1 105&13&12&382 614&1 451 788&525 462\\
6&45 329&245&214&1 547 632&6 132 639&2 462 911\\
\hline
\hline
$\times$&\multicolumn{3}{|c|}{Time (s)}&\multicolumn{3}{|c|}{Memory (B)}\\
\hline\hline
$n$&\texttt{Finotti}&\texttt{phantom}&\texttt{Illusie}&\texttt{Finotti}&\texttt{phantom}&\texttt{Illusie}\\
\hline\hline
1&$<$1&$<$1&$<$1&159 396&502 358&30 987\\
2&$<$1&$<$1&$<$1&172 262&507 208&35 547\\
3&$<$1&$<$1&$<$1&219 313&535 294&60 337\\
4&4&2&2&249 804&819 869&198 237\\
5&207&61&47&385 882&2 444 362&898 503\\
6&11 354&2 108&1 514&1 510 824&11 597 852&4 400 085\\
\hline
\end{tabular}
\caption{Comparisons in $W_n\mleft(\mathbb F_{5^2}\mleft[X\mright]\mright)$ with $d=2$.}\label{n}
\end{figure}

Figure \ref{n} describes the behaviour of the algorithms for varying $n$.
Here the computation time is exponential in $n$.

For the addition, algorithm \ref{illusie} is more than $20$ times faster than Finotti's.
For the multiplication, it is almost $8$ times faster.

\begin{figure}[ht]
\begin{tabular}{|c|c|c|c|c|c|c|}
\hline
$+$&\multicolumn{3}{|c|}{Time (s)}&\multicolumn{3}{|c|}{Memory (B)}\\
\hline\hline
$p$&\texttt{Finotti}&\texttt{phantom}&\texttt{Illusie}&\texttt{Finotti}&\texttt{phantom}&\texttt{Illusie}\\
\hline\hline
2&0&0&0&97 706&530 749&118 882\\
3&0&0&0&125 442&652 118&203 270\\
5&16&20&13&314 978&1 284 997&666 378\\
7&300&219&142&648 530&2 704 277&1 813 790\\
11&9 040&2 333&1 658&3 312 442&8 973 623&6 919 944\\
13&35 063&2 483&1 464&4 248 558&14 518 886&11 841 828\\
\hline
\hline
$\times$&\multicolumn{3}{|c|}{Time (s)}&\multicolumn{3}{|c|}{Memory (B)}\\
\hline\hline
$p$&\texttt{Finotti}&\texttt{phantom}&\texttt{Illusie}&\texttt{Finotti}&\texttt{phantom}&\texttt{Illusie}\\
\hline\hline
2&0&0&0&100 363&572 814&146 310\\
3&0&2&1&124 634&819 853&311 954\\
5&4&43&34&305 440&2 008 946&1 159 036\\
7&25&381&177&857 738&4575 482&3 230 842\\
11&566&5 198&3 414&3 201 702&15 891 410&12 597 546\\
13&2 061&13 524&3 661&4 115 845&25 816 438&21 702 750\\
\hline
\end{tabular}
\caption{Comparisons in $W_4\mleft(\mathbb F_p\mleft[X\mright]\mright)$ with $d=10$.}\label{p}
\end{figure}

Figure \ref{p} describes the behaviour of the three algorithms for varying $p$.
The computation time is also exponential in $p$.

Algorithm \ref{illusie} is more than $200$ times faster than Finotti's for the addition.
For the multiplication however, Finotti's algorithm is almost twice as fast.

\begin{figure}[ht]
\begin{tabular}{|c|c|c|c|c|c|c|}
\hline
$+$&\multicolumn{3}{|c|}{Time (s)}&\multicolumn{3}{|c|}{Memory (B)}\\
\hline\hline
$q$&\texttt{Finotti}&\texttt{phantom}&\texttt{Illusie}&\texttt{Finotti}&\texttt{phantom}&\texttt{Illusie}\\
\hline\hline
$7$&107&46&27&409 210&1 575 463&949 789\\
$7^2$&502&58&40&426 862&1 675 424&721 575\\
$7^3$&787&31&19&427 234&1 675 364&721 457\\
$7^4$&830&50&48&620 470&1 683 176&721 227\\
$7^5$&1 099&55&41&620 634&1 682 817&720 818\\
$7^6$&1 806&63&52&769 907&1 697 348&696 394\\
$7^7$&2 041&64&53&770 172&1 698 681&696 466\\
$7^8$&1 724&33&58&769 996&1 699 734&696 534\\
$7^9$&2 405&71&59&768 482&1 699 615&696 610\\
$7^{10}$&3 580&84&71&770 670&1 699 397&696 746\\
$7^{11}$&2 998&63&55&769 659&1 699 451&696 818\\
$7^{12}$&2 867&65&59&768 865&1 699 414&696 890\\
$7^{13}$&3 228&75&67&769 818&1 698 823&696 962\\
$7^{14}$&2 715&38&38&769 772&1 699 651&697 030\\
$7^{15}$&3 343&54&50&770 064&1 700 071&697 106\\
$7^{16}$&4 268&83&76&769 284&1 699 317&697 174\\
$7^{17}$&4 330&95&87&770 112&1 699 437&697 246\\
$7^{18}$&3 671&86&81&769 397&1 700 182&697 386\\
$7^{19}$&4 555&86&82&768 776&1 699 355&697 458\\
\hline
\hline
$\times$&\multicolumn{3}{|c|}{Time (s)}&\multicolumn{3}{|c|}{Memory (B)}\\
\hline\hline
$q$&\texttt{Finotti}&\texttt{phantom}&\texttt{Illusie}&\texttt{Finotti}&\texttt{phantom}&\texttt{Illusie}\\
\hline\hline
$7$&8&99&62&484 910&2 618 262&1 657 221\\
$7^2$&60&105&69&512 485&2 796 119&1 216 317\\
$7^3$&103&144&93&511 750&2 796 836&1 216 609\\
$7^4$&130&142&100&718 738&2 806 743&1 217 995\\
$7^5$&122&117&103&718 729&2 806 450&1 218 134\\
$7^6$&207&151&110&747 348&2 823 070&1 216 138\\
$7^7$&206&122&96&746 913&2 823 393&1 216 210\\
$7^8$&227&83&56&747 149&2 823 469&1 216 278\\
$7^9$&166&175&126&747 256&2 823 594&1 216 354\\
$7^{10}$&336&191&144&747 743&2 824 107&1 216 490\\
$7^{11}$&243&83&122&746 817&2 823 821&1 216 562\\
$7^{12}$&332&166&125&746 838&2 824 137&1 216 634\\
$7^{13}$&342&188&139&745 731&2 823 978&1 216 706\\
$7^{14}$&366&194&145&747 418&2 824 705&1 216 774\\
$7^{15}$&385&181&139&746 727&2 823 927&1 216 850\\
$7^{16}$&405&202&92&746 825&2 824 757&1 216 918\\
$7^{17}$&543&207&149&747 064&2 823 812&1 216 990\\
$7^{18}$&453&209&118&746 383&2 824 596&1 217 010\\
$7^{19}$&411&166&110&747 615&2 824 175&1 217 082\\
\hline
\end{tabular}
\caption{Comparisons in $W_4\mleft(\mathbb F_q\mleft[X\mright]\mright)$ with $d=5$ and $p=7$.}\label{q}
\end{figure}

In figure \ref{q} we compare how greater powers of $p$ for the cardinality of the base field have an impact on the algorithms.
We see that, especially for the addition, Finotti's algorithm is more affected.

This is something which also happens for multivariate polynomials even though, in that situation, Finotti's approach is very efficient.

\begin{figure}[ht]
\begin{tabular}{|c|c|c|c|c|c|c|}
\hline
$+$&\multicolumn{3}{|c|}{Time (s)}&\multicolumn{3}{|c|}{Memory (B)}\\
\hline\hline
$q$&\texttt{Finotti}&\texttt{phantom}&\texttt{Illusie}&\texttt{Finotti}&\texttt{phantom}&\texttt{Illusie}\\
\hline\hline
$5$&70&1 013&962&3 950 636&32 996 786&17 968 328\\
$5^2$&1 157&1 076&1 213&3 158 069&36 165 400&19 727 049\\
\hline
\hline
$\times$&\multicolumn{3}{|c|}{Time (s)}&\multicolumn{3}{|c|}{Memory (B)}\\
\hline\hline
$q$&\texttt{Finotti}&\texttt{phantom}&\texttt{Illusie}&\texttt{Finotti}&\texttt{phantom}&\texttt{Illusie}\\
\hline\hline
$5$&12&10 267&5 958&7 017 630&118 633 580&57 119 577\\
$5^2$&306&17 516&4 241&5 815 799&135 236 744&64 295 772\\
\hline
\end{tabular}
\caption{Comparisons in $W_4\mleft(\mathbb F_q\mleft[X,Y\mright]\mright)$ with total degrees of the coefficients at most $2$.}\label{qq}
\end{figure}

One can see in figure \ref{qq} how quickly, especially for the addition, Finotti's algorithm is affected by the growth of $q$.
However, Finotti's approach is so efficient in the multivariate case, especially for the multiplication, that in many cases it outperforms both new methods presented in this paper.

\bigskip

To sum up, there are five parameters which will affect the computation of Witt vectors : the maximum (total) degree $d$ of the coefficients, the length $n$ of the truncated Witt vectors, the characteristic $p$ and cardinality $q$ of the base field, and the number $x$ of indeterminates of the polynomial ring.

For memory usage, Finotti's method is almost always better.

For speed, the parameters $d$ and $n$ do not seem to affect much the choice of the algorithm.
For bigger $q$, and to a lesser extent for bigger $p$, algorithm \ref{illusie} is to be preferred, and one should also use it for small $x$ (especially for $x=1$).

Algorithm \ref{illusie} is much faster for addition, whereas for multiplication it depends on the parameters.

\end{document}